\documentclass[12pt,letterpaper]{amsart}
\usepackage{amssymb,latexsym, amsmath, amsxtra}

\theoremstyle{plain}

\theoremstyle{definition}

\newtheorem*{exercise}{Exercise}
\theoremstyle{remark}

\numberwithin{equation}{section}

\def\({\left(}
\def\){\right)}

\def\yr #1#2#3#4{(#1#2#3#4)}
\def\publ{\relax}
\def\publaddr{\relax}
\def\vol{\relax}

\def\jour{\relax}
\def\pages{\relax}

\begin{document}
\title{Basic analytic number theory}
\author{David W Farmer}

\thanks{Work supported by the
American Institute of Mathematics
and by the Focused Research Group grant (0244660) from the NSF.  
This paper will appear in the proceedings of the school
``Recent Perspectives in Random Matrix Theory and Number Theory''
held at the
Isaac Newton Institute, April 2004.
}

\thispagestyle{empty}
\vspace{.5cm}
\begin{abstract}
We give an informal introduction to the most basic techniques
used to evaluate moments on the critical line of the
Riemann zeta-function and to find asymptotics for sums of
arithmetic functions.  
\end{abstract}

\address{
{\parskip 0pt
American Institute of Mathematics\endgraf
360 Portage Ave.\endgraf
Palo Alto, CA 94306\endgraf
farmer@aimath.org\endgraf
}
  }

\maketitle

\section{Introduction}

The simplest way to compute a moment of the zeta-function is to 
approximate the zeta-function by Dirichlet polynomials and then
compute the moment of the polynomials.  In this paper we 
describe the most rudimentary techniques in this area.  Along the way
discuss the basic methods for finding asymptotics of sums of arithmetic
functions, and we also compute the arithmetic factor in the standard
conjectures for moments of the zeta-function.  The methods described here
are completely standard: our intention is to give a brief
introduction to those who are new to the subject.
The standard reference is
Titchmarsh~\cite{T} and we cite the specific sections where one can look 
for more details.

We assume a knowledge of the calculus of one complex variable.
Readers unfamiliar with the big-$O$, little-$o$, and $\ll$ notation
(and physicists, who may use the $\ll$ notation differently),
should consult the Appendix.

\section{The ``first moment''}

The simplest approximation to the zeta-function is
\begin{equation}
\zeta(s)=\sum_{1\le n\le T} \frac{1}{n^s} - \frac{T^{1-s}}{1-s} 
+ O(T^{-\sigma}),
\end{equation}
where $s=\sigma+it$ and the equation is valid for $|t|\le T$.
See~\cite{T}, Section 4.11. 
Specializing to $s=\frac12+it$, we have
\begin{equation}\label{eqn:sumapproximation}
\zeta(\tfrac12+it)=\sum_{1\le n\le T} \frac{1}{n^{\frac12+it}} 
+O\(\frac{T^{\frac12}}{1+|t|}\).
\end{equation}
Now we can find the average value of the zeta-function on the
critical line.  The justification of the steps follows the calculation.
\begin{eqnarray}
\int_0^T\zeta(\tfrac12+it) \,dt&=&
\sum_{1\le n\le T} \frac{1}{\sqrt{n}} \int_0^T n^{-it}dt 
+ O\(T^{\frac12}\int_0^T \frac{1}{1+|t|}\, dt\) \cr
&=&T + \sum_{2\le n\le T}\frac{1}{\sqrt{n}}\frac{1-n^{-iT}}{\log n} 
+ O\(T^{\frac12}\log T\) \cr
&=& T+ O\(T^{\frac12}\log T\) .
\end{eqnarray}
In the first step we switched the sum and the integral, which is justified 
because both are finite.  The sum on the second line was estimated by
\begin{equation}
\sum_{2\le n\le T}\frac{1}{\sqrt{n}}\frac{1-n^{-iT}}{\log n}
\ll \sum_{2\le n\le T} \frac{1}{\sqrt{n}\log n}
\ll \int_2^T \frac{1}{\sqrt{x}\log x} dx
\ll \frac{T^{\frac12}}{\log T},
\end{equation}
which is smaller than the other error term.  

Thus, the zeta-function on the critical line is 1 on average.  This is
not a particularly useful piece of information, for it is the magnitude
of the zeta-function which is of interest.  So we consider the
moments of $|\zeta(\frac12+it)|$.

\section{The 2nd moment}

The simple fact that $|z|^2=z\overline{z}$ means that 
$|\zeta(\frac12+it)|^K$ is much more amenable to methods
of complex analysis when $K=2k$ is an even integer.  
The easiest case is the 2nd moment, which we compute in this section.

For later use it will be helpful to first consider the 2nd moment of a general
Dirichlet polynomial.
Suppose
\begin{equation}
P(s)=\sum_{1\le n\le N}\frac{a_n}{n^{s}} .
\end{equation}
We have
\begin{eqnarray}
\int_0^T \left|P(it)\right|^2 \, dt
&=&\int_0^T \left|\sum_{1\le n\le N}\frac{a_n}{n^{it}}\right|^2 dt \cr
&=& \int_0^T \sum_n \frac{a_n}{n^{it}} \sum_m \frac{\overline{a_m}}{m^{-it}} \, dt \cr
&=& \sum_{n,m} a_n \overline{a_m} \int_0^T \(\frac{m}{n}\)^{it} \, dt \cr
&=& T \sum_n |a_n|^2 
    +  \sum_{n\not =m} \frac{a_n \overline{a_m} \(\(\frac{m}{n}\)^{iT}-1\) }
					{\log(m/n)} \cr
&=& T\, \mathcal M(N)+\mathcal E(N,T),
\end{eqnarray}
say.  We think of $T\, \mathcal M(N)$ as the main term and $\mathcal E(N,T)$ as an error term,
so we want to understand when $\mathcal E(N,T)$ will be smaller in magnitude than
$T\, \mathcal M(N)$.  

Setting $m=n+h$ we can rewrite the error term as
\begin{equation}
\mathcal E(N,T)=
\sum_n\sum_{h\not = 0}  
\frac{a_n \overline{a_{n+h}} \(\(1+\frac{h}{n}\)^{iT}-1\) }
                                        {\log(1+h/n)} .
\end{equation}
Now consider only the $h=1$ term from the above sum:
\begin{equation}
\sum_{1\le n\le N} a_n \overline{a_{n+1}} \cdot n \cdot (\text{something bounded}) .
\end{equation}
Without any information on $a_n$, the above sum, which is just one part of
$\mathcal E(N,T)$, could be about the same size as~$N\, \mathcal M(N)$.  
Thus, in general one should only expect
$\mathcal E(N,T)$ to be smaller than $T\, \mathcal M(N)$ if $N<T$.  And if $N>T$ then one 
may need some detailed
information about the coefficients~$a_n$ in order 
to extract something meaningful from~$\mathcal E(N,T)$. Goldston and Gonek~\cite{GG}
have given a clear discussion of these issues.

One can carry through the above calculation to obtain a useful general
mean value theorem
for Dirichlet polynomials.  See Titchmarsh~\cite{T}, Section~7.20.  However,
a stronger result is provided 
by the mean value theorem of Montgomery and Vaughan~\cite{MV}:
\begin{equation}\label{eqn:mv}
\int_0^T \left|\sum_{1\le n\le N}\frac{a_n}{n^{it}}\right|^2 dt
=
\sum_{1\le n\le N} |a_n|^2 \(T+O(n)\) .
\end{equation}
This result is best possible for general sequences~$a_n$. 
Note that if $N=o(T)$ then the error term is smaller than the main term.

To use the mean value theorem to compute the second moment
of the zeta-function, first write
(\ref{eqn:sumapproximation}) as $\zeta=S+E$.  
That is,
\begin{equation}
S=\sum_{1\le n\le T} \frac{1}{n^{\frac12+it}} 
\ \ \ \ \ \ \ \ \ \ \ \ 
\text{and}
\ \ \ \ \ \ \ \ \ \ \ \ 
E=O\(\frac{T^{\frac12}}{1+|t|}\).
\end{equation}
Using 
\begin{equation}
|\zeta|^2=|S+E|^2 = 
(S+E)(\overline{S}+\overline{E})=|S|^2+2\Re S\overline{E} + |E|^2,
\end{equation}
we have
\begin{equation}
\int_0^T\left| \zeta(\tfrac12+it)\right|^2dt = 
\int_0^T\left|S\right|^2dt +
2 \Re \int_0^T S\overline{E} \, dt
+
\int_0^T\left|E\right|^2dt .
\end{equation}
We want to evaluate the $|S|^2$ integral as our main term and estimate
the $|E|^2$ integral as our error term, but what to do about the 
cross term?  By the Cauchy-Schwartz inequality,
\begin{equation}
2 \Re  \int_0^T S\overline{E} \, dt \ll  \int_0^T |S||E| \, dt
\ll \(\int_0^T |S|^2 dt\)^\frac12 \(\int_0^T |E|^2 dt\)^\frac12 .
\end{equation}
If $\int_0^T |E|^2 dt$ is smaller than $\int_0^T |S|^2 dt$, then so
is the right side of the inequality above.  We have the general
principle that if the ``main error term'' is smaller than the main
term, then so are the cross terms.  It remains only to evaluate
$\int_0^T |S|^2 dt$ and estimate $\int_0^T |E|^2 dt$.

For the main term use~(\ref{eqn:mv}) with $a_n=n^{-\frac12}$:
\begin{eqnarray}
\int_0^T |S|^2 dt &=&\sum_{1\le n\le T} \frac1n \(T+O(n)\) \cr
&=& T\(\log T + O(1)\) + \sum_{1\le n\le T} O(1) \cr
&=& T\log T + O(T).
\end{eqnarray}
For the error term we have
\begin{equation}
\int_0^T |E|^2 dt \ll T\int_0^T \frac{1}{(1+|t|)^2} \, dt \ll T,
\end{equation}
which is smaller than the main term.  We have just proven the mean value
result
\begin{equation}\label{eqn:second}
\int_0^T|\zeta(\tfrac12+it)|^2dt = T\log T + O(T\log^{\frac12} T).
\end{equation}
Note that we do not obtain an error term 
of $O(T)$ in (\ref{eqn:second}) 
by these methods.  However, a much better error term can be obtained.
The first step toward this is given in the next section.

\begin{exercise}  Approximate the sum by an integral to show that
if $A>1$ then
\begin{equation}
\sum_{n>T} \frac{1}{n^A} = \frac{T^{1-A}}{A-1} + O(T^{-A}).
\end{equation}
Conclude that if $A>1$ then
\begin{equation}
\sum_{n\le T} \frac{1}{n^A} = \zeta(A)-\frac{T^{1-A}}{A-1} + O(T^{-A}).
\end{equation}
\end{exercise}

\begin{exercise} Show that if $\frac12<\sigma<1$ then 
\begin{equation}\label{eqn:secondsigma}
\int_0^T|\zeta(\sigma+it)|^2dt = T\(\zeta(2\sigma)-\frac{T^{1-2\sigma}}{2\sigma-1}\)
 + O\(\frac{T^{\frac32 -\sigma}}{1-\sigma} \log^{\frac12} T\),
\end{equation}
where the implied constant in the big-$O$ term is independent of~$\sigma$.
\end{exercise}

Since the error term in (\ref{eqn:secondsigma}) is uniform in sigma, we
can let $\sigma\to\frac{1}{2}^+$
to recover~(\ref{eqn:second}).  This makes use of the fact that 
$\zeta(s)$ has a simple pole with residue~1 at~$s=1$. 

Note that if $\sigma>\frac12$ is independent of $T$ then the right side
of (\ref{eqn:secondsigma}) is of size $\approx T$.  On the other hand, 
the second moment on the $\frac12$-line is of size $T\log T$.  
Thus there is an abrupt change in the behavior of the zeta-function 
when one moves onto the critical line.  Equation (\ref{eqn:secondsigma})
illustrates that the transition occurs on the scale of
$1/\log T$ from the $\frac12$-line.

\section{Better 2nd moment}

The methods of the previous section are not sufficient to evaluate
the main term of the 2nd moment with an error term $O(T^A)$ for
$A<1$, nor are those methods sufficient to evaluate
the 4th moment of the zeta-function.  
Evaluating the 4th moment by squaring (\ref{eqn:sumapproximation})
gives a Dirichlet polynomial of length $T^2$, which cannot be handled
by the Montgomery-Vaughan mean value theorem.  So one needs either
a shorter approximation to $\zeta(s)$, or an approximation to
$\zeta^2(s)$ of length $\le T$, or a way to handle longer polynomials.
In preparation for the 4th moment, we first evaluate the
2nd moment with a better error term.

The ``approximate functional equation'' of Hardy and Littlewood
expresses the $\zeta$-function as a sum of \emph{two} short
Dirichlet polynomials:
\begin{equation}
\zeta(s)=\sum_{1\le n\le x} \frac{1}{n^s} 
      + \chi(s) \sum_{1\le n\le y} \frac{1}{n^{1-s}}
+ O\(x^{-\sigma}+|t|^{\frac12-\sigma}y^{\sigma-1}\),
\end{equation}
where $xy=t/2\pi$ and $\chi(s)$ is the usual factor in the
functional equation $\zeta(s)=\chi(s)\zeta(1-s)$.  The
name ``approximate functional equation'' comes from the fact
that the right side looks like $\zeta(s)$ when $x$ is large
and like $\chi(s)\zeta(1-s)$ when $y$ is large.
On the $\frac12$-line we have
\begin{equation}\label{eqn:approxhalf}
\zeta(s)=\sum_{1\le n\le N} \frac{1}{n^{\frac12+it}} 
      + \chi(s) \sum_{1\le n\le N} \frac{1}{n^{\frac12-it}}
+ O\(N^{-\frac12}\),
\end{equation}
where we set $x=y=N=\sqrt{t/2\pi}$.

Hardy and Littlewood used the approximate functional equation
to evaluate the second moment of the zeta-function with a better
error term.
We will not carry out the calculation, but just give
the flavor.  See Titchmarsh~\cite{T} Section~7.4 for details.
Writing (\ref{eqn:approxhalf}) as
$\zeta=S+\chi(s)\overline{S}+E$,
we have
\begin{equation}
\int_0^T |\zeta(\tfrac12+it)|^2 dt =
2\int_0^T |S|^2 dt + 2 \Re \int_0^T \chi(\tfrac12-it) S^2 dt
+ \text{error term}.
\end{equation}
In that calculation we used the fact that $|\chi(\tfrac12+it)|=1$.
Below we will use the more precise information
\begin{equation}\label{eqn:chiequals}
\chi(s)=\(\frac{t}{2\pi}\)^{\frac12-s}e^{i t+\pi i/4} \(1+O(t^{-1})\).
\end{equation}

The main difficulty is evaluating the $\chi(\tfrac12-it) S^2$ term,
which equals
\begin{equation}
\sum_{n,m}\frac{1}{\sqrt{nm}} \int_0^T
\chi(\tfrac12-it) (nm)^{-it} dt.
\end{equation}
By (\ref{eqn:chiequals}), up to a negligible error that integral is of the form 
$\int e^{i f(t)} dt$.  If $f(t)$ has a stationary point in the
range of integration then we can extract a main term, otherwise
it will become an error term. 
In particular, that integral can be handled by the method
of stationary phase.
 See Titchmarsh~\cite{T} Section~7.4
for details.  Our point here is that by virtue of
(\ref{eqn:chiequals}), integrals involving $\chi(\frac12+it)$ and
Dirichlet polynomials can be handled.
The result that can be obtained by the above argument is
\begin{equation}
\int_0^T |\zeta(\tfrac12+it)|^2 dt = T\log T + (2\gamma-1)T 
+ O(T^{\frac12+\varepsilon}).
\end{equation}
The error term can be improved by more sophisticated methods.

There are many applications of mean values of the zeta-function
multiplied by a Dirichlet polynomial.  
Suppose
\begin{equation}
M(s)=\sum_{1\le n <T^\theta} \frac{b_n}{n^s} 
\end{equation}
with $b_n\ll n^\varepsilon$.
By the approximate functional equation, $\zeta(s)M(s)$ can be approximated
by Dirichlet polynomials of length~$T^{\frac12+\theta}$.
So by the above methods, if $\theta<\frac12$ then one should be able to find
an asymptotic formula for $\int_0^T |\zeta M(\tfrac12+it)|^2 dt$.
Evaluating such an integral, with $\theta=\frac12-\varepsilon$, 
was key to Levinson's proof that more 
than one-third of the zeros of the zeta-function are on the critical line.
Conrey made use of very deep and technical results to evaluate such an
integral with $\theta=\frac47-\varepsilon$, leading to the result that
more than two-fifths of the zeros are on the critical line.

\section{The 4th moment}
To evaluate the 4th moment of the zeta-function by the methods described above, 
one requires an approximation to $\zeta(s)^2$ of length less than~$T$.
This is provided by the following approximate functional equation:
\begin{equation}\label{eqn:zetasquaredapfe}
\zeta(s)^2=\sum_{1\le n\le x}\frac{d(n)}{n^s} 
   + \chi(s)^2 \sum_{1\le n\le y} \frac{d(n)}{n^{1-s}}
+ O(x^{\frac12-\sigma}\log t), 
\end{equation}
where $x y=(t/2\pi)^2$ and $d(n)$ is the number of divisors of~$n$.
More generally one should expect an approximate functional equation
of the form
\begin{equation}\label{eqn:zetakapfe}
\zeta(s)^k=\sum_{1\le n\le x}\frac{d_k(n)}{n^s} 
   + \chi(s)^k \sum_{1\le n\le y} \frac{d_k(n)}{n^{1-s}}
+ \text{ error term}
\end{equation}
where $x y\approx t^{k}$ and $d_k(n)$ is the $k$-fold divisor
function
\begin{equation}
d_k(n) = \sum_{n_1\cdots n_k=n} 1,
\end{equation}
which has generating function
\begin{equation} \label{eqn:zetak}
\zeta(s)^k=\sum_{n=1}^\infty \frac{d_k(n)}{n^s},
\ \ \ \ \ \ \ \ \sigma>1.
\end{equation}

Plugging (\ref{eqn:zetasquaredapfe}) into the Montgomery-Vaughan
mean value theorem leads to
\begin{eqnarray}
\int_0^T |\zeta(\tfrac12+it)|^4 dt = \frac{1}{2\pi^2}T\log^4 T +
\text{ error term}.
\end{eqnarray}
If you actually do the calculation, you will find that in order
to determine the main term you need to evaluate
sums like
\begin{equation}\label{eqn:dsquaredsum}
\sum_{1\le n\le X} \frac{d(n)^2}{n} .
\end{equation}
There is a standard technique for finding the leading-order
asymptotics of such sums, which is given in the next section.

\subsection{Comments on approximate functional equations}

The error term in (\ref{eqn:zetasquaredapfe}) is rather large
and leads to an error term of size $O(T\log^2 T)$ in
the 4th moment $\int_0^T |\zeta(\frac12+it)|^4 dt$. 
The large size of the error term is due to the fact that our sums have a
sharp cut-off. The error term can
be much reduced by having a smooth weight in the sums.
That is, 
\begin{equation}
\zeta(s)^k=\sum_{n}\frac{d_k(n)}{n^s} \varphi(n,t)
   + \chi(s)^k \sum_{n} \frac{d_k(n)}{n^{1-s}}\varphi^*(n,t)
+ \text{ small error term},
\end{equation}
where $\varphi$ and $\varphi^*$ are particular functions 
that are approximately 1 for $n<t^{\frac{k}{2}}$ and decay
for $n>t^{\frac{k}{2}}$.  
As our previous discussion should suggest, it is the 
length of the sums, and not the size of the error term,
which provides the true difficulty when~$k>2$.

\section{Perron's formula}

Here is the problem: you have an arithmetical function $a_n$ and
you want to find the asymptotics of
\begin{equation}\label{eqn:ansum}
S(X)=\sum_{1\le n\le X} a_n .
\end{equation}
This problem can often be solved by the most basic methods of
analytic number theory.

First note that for integers $N\ge 1$,
\begin{equation}\label{eqn:basicperron}
\frac{1}{2\pi i}\int_{1-i Y}^{1+ i Y} A^s\, \frac{ds}{s^N} =
\begin{cases}
\frac{\log^{N-1} A}{(N-1)!}  + \text{ error term} & \text{ if } A>1 \\
0 + \text{ error term} & \text{ if } A<1
\end{cases}
\end{equation}
To see this, consider the integral
\begin{equation}
\frac{1}{2\pi i}
\int_{{\mathcal C}_1} A^s \,\frac{ds}{s^N} 
\end{equation}
where the integration is over the closed rectangular path connecting the points
\begin{equation}
\mathcal C_1 = 
\begin{cases}
[1-i Y, 1+ i Y,-B+i Y, -B-iY], & \text{ if } A>1 \\
[1-i Y, 1+ i Y,B+i Y, B-iY], &  \text{ if } A<1 ,
\end{cases}
\end{equation}
where $B$ is a large positive number.
In both cases the main term comes from the residue of the pole at~$0$,
which is or is not inside the path of integration.  

\begin{exercise} Bound the error term in (\ref{eqn:basicperron}) by estimating
the integral along the three segments of $\mathcal C_1$ other than
$[1-i Y, 1+ i Y]$. You should find that if $N\ge 2$,
 then you can let $Y\to\infty$
and the error vanishes.
See Section 3.12 of \cite{T} if you aren't sure how to begin.
\end{exercise}

To evaluate (\ref{eqn:ansum}), let
\begin{equation}
F(s)=\sum_{n=1}^\infty \frac{a_n}{n^s},
\end{equation}
and suppose that the sum converges absolutely for $\sigma>\sigma_0$.
Using (\ref{eqn:basicperron}) and supposing $\sigma>\sigma_0$, we have
\begin{eqnarray}\label{eqn:perron}
\frac{1}{2\pi i}
\int_{\sigma-i Y}^{\sigma+ i Y} F(s) X^s \, \frac{ds}{s} 
&=&
\sum_{n=1}^\infty a_n \,
\frac{1}{2\pi i}
\int_{\sigma-i Y}^{\sigma+ i Y} \(\frac{X}{n}\)^s \frac{ds}{s} \cr
&=&
\sum_{n=1}^\infty a_n 
\begin{cases} 1+ \text{ error term} & \text{ if } X>n \\
0 + \text{ error term} & \text{ if } X<n
\end{cases}
\cr
&=& S(X) + \text{ error term}.
\end{eqnarray}
This is known as Perron's formula.
If we can learn enough about the function $F(s)$ so that
the integral in (\ref{eqn:perron}) can be evaluated in another way,
then we will have a formula for~$S(X)$.

Suppose (as is frequently the case) that $F(s)$ has a pole at
$\sigma_0$ and no other poles in the
half-plane $\sigma>\sigma_1$ for some $\sigma_1<\sigma_0$.
Then consider
\begin{equation} \label{eqn:contourintegral}
\frac{1}{2\pi i}
\int_{{\mathcal C}_2(\varepsilon)}
F(s) X^s \, \frac{ds}{s},
\end{equation}
where for $\varepsilon>0$ the integration is over the rectangular path with vertices
\begin{equation}
\mathcal C_2(\varepsilon) = 
[\sigma-i Y, \sigma+ i Y,
\sigma_1+\varepsilon+i Y, \sigma_1+\varepsilon-iY].
\end{equation}
We can evaluate (\ref{eqn:contourintegral}) by finding the residue at the pole
$s=\sigma_0$ and at $s=0$ (if 0 is inside the path of integration.
And in the same way as you estimated the error term in (\ref{eqn:basicperron}),
we find that (\ref{eqn:contourintegral}) equals the integral
in Perron's formula plus an error term. The final step of bounding the 
integral on the 3 other segments requires the additional ingredient of a 
bound for $F(\sigma+it)$ as $t\to\infty$, uniformly for $\sigma>\sigma_1$.

For example, at the end of the previous section we wanted to evaluate
the sum of $d(n)^2/n$.  

\begin{exercise} Check that
\begin{equation}
\sum_{n} \frac{d(n)^2}{n^s} = \frac{\zeta^4(s)}{\zeta(2s)}.
\end{equation}
Hint: both sides have an Euler product.  The factors on the right can be 
found from the Euler product for the zeta-function.  Those on the left
require summing $\sum_{j=0}^\infty d(p^j)^2 p^{-js}$.
\end{exercise}

Thus, we apply Perron's formula~(\ref{eqn:perron}) with
\begin{equation}
F(s) = \frac{\zeta^4(s+1)}{\zeta(2s+2)}.
\end{equation}
To determine the analytic properties of $F(s)$, use the fact that $\zeta(s)$
is entire except for a simple pole at $s=1$, where we have the Laurent expansion
\begin{equation}
\zeta(s) = \frac{1}{s-1} + \gamma + \cdots .
\end{equation}
Also $\zeta(s)$ has no zeros in $\sigma>1$, and no zeros in  $\sigma>\frac12$
assuming the Riemann Hypothesis.  In these calculations one frequently 
needs that $\zeta(2)=\pi^2/6$  and
\begin{equation}
\zeta(s)= -\tfrac12 -\tfrac12 \log(2\pi)s  + \cdots  
\end{equation}
for $s$ near~$0$.

To estimate the error terms,  one can use  the convexity estimate
\begin{equation}
\zeta(\sigma+it)
\ll
\begin{cases}
1 & \sigma>1 \\
|t|^{\frac12 - \frac12 \sigma} & 0<\sigma<1 \\
|t|^{\frac12 - \sigma} & \sigma< 0, \\
\end{cases}
\end{equation}
along with $\zeta(s)\gg1$ for $\sigma>1$. 
Also, assuming RH we have
$t^{-\varepsilon}\ll \zeta(\sigma+it) \ll t^{\varepsilon}$ for 
$\sigma>\frac12$. All of these estimates are for fixed~$\sigma$
as $t\to\infty$.

Assembling the pieces we find
\begin{equation}\label{eqn:dn2answer}
\sum_{1\le n\le X} \frac{d(n)^2}{n} \sim \frac{\log^4 X}{4\pi^2}.
\end{equation}

\begin{exercise}Argue that (\ref{eqn:dn2answer}) is of
the form $X P_4(\log X)+O(X^B)$ where $P_4$ is a polynomial of
degree~4 and $B<1$. Find the next-to-leading coefficient of $P_4$, and estimate 
$B$ both with and without assuming the Riemann Hypothesis.
\end{exercise}

\begin{exercise}Deduce the following asymptotics:
\begin{eqnarray}
\sum_{1\le n\le X} d_k(n) &\sim& \frac{1}{k!} X\log^k X \cr
\sum_{1\le n\le X} \varphi(n) &\sim& \frac{3}{\pi^2} X^2,
\end{eqnarray}
where $\varphi(n)$ is the Euler totient function.
In addition to the generating
function (\ref{eqn:zetak}), you should use (and prove):
\begin{equation}
\sum_{n} \frac{\varphi(n)}{n^s} = \frac{\zeta(s-1)}{\zeta(s)}.
\end{equation}
\end{exercise}

\begin{exercise} Find the next-to-leading order terms in the previous
exercise. Also determine the shape of the main terms and estimate the
size of the error terms, both with and without the Riemann Hypothesis.
\end{exercise}

Note that there are interesting and important cases where the above analysis
is inadequate.  For example, in the proof of the prime number theorem 
$a_n=\Lambda(n)$, the von~Mangoldt function.
Then $F(s)$ has a pole at $s=1$ as well as
poles at the zeros of the $\zeta$-function and one must
use a more complicated path of integration as well as nontrivial 
estimates for $\zeta(s)$ in the critical strip.

\section{The conjecture for moments of the zeta-function}

Much recent work on the relationship between $L$-functions and Random Matrix Theory
was motivated by the problem of finding conjectures for the $2k$th moment of
the Riemann zeta-function on the critical line.  Conrey and Ghosh~\cite{CG}
formulated it as follows: for each integer $k\ge0$ there exists an integer $g_k$
such that 
\begin{equation}
\int_0^T |\zeta(\tfrac12+it)|^{2 k} dt \sim
g_k\, 
\frac{a_k}{k^2 !} \,T\log^{k^2} T,
\end{equation}
where
\begin{equation}
a_k=\prod_p \(1-\frac1p\)^{k^2} \sum_{m=0}^\infty \binom{k+m-1}{m} p^{-m} .
\end{equation}
In this conjecture the only missing ingredient is the \emph{integer}~$g_k$.
Keating and Snaith~\cite{KS} computed the moments of characteristic polynomials
of unitary matrices and used the result to conjecture
\begin{equation}
g_k=k^2 ! \prod_{j=0}^{k-1} \frac{j!}{(k+j)!} .
\end{equation}
It is not trivial to show that this $g_k$ is actually an integer~\cite{CF}.

Our last topic in this paper is to show how the factor $a_k$ arises naturally.
From the approximate functional equation (\ref{eqn:zetakapfe}) it is reasonable
to consider
\begin{equation}
\int_0^T \left|\sum_{n<T} \frac{d_k(n)}{n^{\frac12+it}}\right|^2 dt.
\end{equation}
Note that the sum has length~$T$.  This is good because we can use the
mean value theorem.  But it is bad because the polynomial is not long 
enough to fully approximate $\zeta(\frac12+it)^k$.  We cannot expect
this mean value to equal the $2k$th moment of the zeta-function,
but how far will it be off?  It would be nice if it were off by
some simple factor, so
one possible interpretation of $g_k$ is
``the number of length $T$ polynomials needed to capture
the $2k$th moment of the $\zeta$-function.''  
We do not claim that this was the original reasoning of Conrey and Ghosh.

By the Montgomery-Vaughan mean value theorem the above integral has main term
\begin{equation}
\sum_{n<T} \frac{d_k(n)^2}{n} .
\end{equation}
By Perron's formula, to evaluate this we need to find the leading pole of
\begin{equation}
F(s)=\sum_{n<T} \frac{d_k(n)^2}{n^s} .
\end{equation}
If $k>2$ then $F(s)$ is not a simple expression involving
known functions, but fortunately we do not require complete information
about~$F(s)$.  

First note that
\begin{eqnarray}
\zeta(s)&=&\prod_p \(1+\frac{1}{p^s} + \cdots \) \cr
&=& \frac{1}{s-1} + \cdots ,
\end{eqnarray}
so
\begin{eqnarray}
\zeta(s)^N&=&\prod_p \(1+\frac{N}{p^s} + \cdots \) \cr
&=& \frac{1}{(s-1)^N} + \cdots .
\end{eqnarray}
Thus, if the coefficients of $p^{-js}$ in an Euler product are
integers that only depend on $j$,
then the coefficient of $p^{-s}$ tells you the
order of the pole at $s=1$.
Since
\begin{eqnarray}
d_k(p)&=&\sum_{n_1\cdots n_k=p}1 \cr
&=& k,
\end{eqnarray}
we see that  $F(s)$ has a pole of order $k^2$ at $s=1$, that is, 
$F(s)=a_k (s-1)^{-k^2}+\cdots$, where~we will show that $a_k$ is as given above.
To see that $F(s)$
has no other poles in $\sigma>\frac12$, note that
\begin{equation}\label{eqn:firstestermann}
\zeta^{-k^2}(s) F(s) = \prod_p \(1+ \frac{\beta_2}{p^{2s}} + \frac{\beta_3}{p^{3s}} + \cdots\),
\end{equation}
where the $\beta_j$ are certain integers that do not grow too fast.
In particular, the above Euler product converges absolutely
for $\sigma>\frac12$ so it represents a regular function that is bounded
in $\sigma>\frac12+\varepsilon$.
We have all of the pieces to apply the methods of the previous section,
giving
\begin{equation}
\int_0^T \left|\sum_{n<T} \frac{d_k(n)}{n^{\frac12+it}}\right|^2 dt
\sim
\frac{a_k}{k^2!}\, T\log^{k^2}T,
\end{equation}
where
\begin{eqnarray}
a_k&=& \lim_{s\to 1} (s-1)^{k^2} F(s)\cr
&=& \lim_{s\to 1} \zeta(s)^{-k^2} F(s) \cr
&=& \prod_p \(1-\frac{1}{p}\)^{k^2} \sum_{m=0}^\infty d_k(p^m)^2 p^{-m} .
\end{eqnarray}
Note that the product converges because $d_k(p)=k$ and $d_k(n)\ll n^\varepsilon$.

Finally, $d_k(p^m)=\binom{k+m-1}{m}$, as can be seen by the following argument.
Since $d_k(p^m)$ is the number of ways of writing $e_1+\cdots+e_k=m$, we can select
the $e_j$ by writing down $m+k-1$ circles~$\circ$ and filling in $k-1$ of them to make
a dot~$\bullet$.  Then
the $e_j$ are the number of circles between the dots, 
including the circles before the 
first and after the last dot.  For example, here is one configuration that
arises from $d_5(p^3)$:
\begin{eqnarray}
&&\mathstrut\hbox{$\ \ \bullet \circ \circ \bullet \circ \bullet \bullet$} \cr
p^3=&&\mathstrut\hbox{$1\ \, \ \ p^2\ \ \ p\ \ 1 \ \ 1 \ $}
\end{eqnarray}

\section{The Estermann phenomenon}

The idea behind formula (\ref{eqn:firstestermann}) can be
generalized to show that if $c(n)$ is a multiplicative
function such that $c(n)\ll n^\varepsilon$ and $c(p^j)$ is an integer that is
independent of the prime~$p$, then
\begin{eqnarray}
F(s)&=&\sum_{n=1}^\infty \frac{c(n)}{n^s} \cr
&=& Z_J(s) \prod_{j<J}  \zeta(js)^{C(j)},
\end{eqnarray}
where
the $C(j)$ are integers and $Z(s)$ is
regular and bounded in $\sigma>1/J$.
Thus $F(s)$, which is originally defined
for $\sigma>1$, has a meromorphic continuation to~$\sigma>0$.

Note that $F(s)$ cannot be continued past $\sigma=0$ 
unless $C(j)=0$ for almost
all $j$. This is because
the zeros of the zeta-function lead to zeros or poles of $F(s)$ 
that accumulate along the $\sigma=0$ line,
giving a natural boundary.  This is known as ``the Estermann phenomenon''.

\section{Appendix: big-$O$ and $\ll$ notation}

The statement
\begin{equation}
f(x)=O(g(x))
\ \ \ \ \ \ \ \ \ \ \ \ 
\text{as }\ \ x\to\infty
\end{equation}
is pronounced ``$f(x)$ is big oh of $g(x)$.''  It is equivalent to
\begin{equation}
f(x) \ll g(x)
\ \ \ \ \ \ \ \ \ \ \ \ \ \ \ \ \ 
\text{as }\ \ x\to\infty,
\end{equation}
which is pronounced
``$f(x)$ is less than less than $g(x)$.''  The symbol $\ll$ is 
typed as \verb1\ll1
in \TeX.  Both of the above statements mean the following: there exists
a constant $C$ such that if $x$ is sufficiently large then
$|f(x)|\le C \, g(x)$.  The number $C$ is called ``the implied constant.''

Note:
\begin{itemize}
\item $f(x) \ll g(x)$ does \emph{not} mean that $f(x)$ is much smaller than
$g(x)$.  It is more accurate to say that $f(x)$ does not grow faster
than $g(x)$.

\item the above statements have the condition ``as $x\to\infty$''.
It is also common to use the big-$O$ and $\ll$ notation to describe
the behavior of a function as $x\to 0$.  Then the definition is modified
to ``if $x$ is sufficiently small''.   Usually context makes it clear
which behavior is being considered.  

\item Both notations are useful: the $\ll$ does not require parentheses,
and the big-$O$ can be used as one term in a formula.
\end{itemize}

Here are some examples.
Below, $A$ and $\varepsilon$
are arbitrary fixed positive numbers.  

Examples assuming $x\to\infty$:
\begin{eqnarray}
x^3 &\ll& x^4 \cr
\log(x)&=& O(x)   \cr
\log(x)&\ll& x^\varepsilon \cr
x^A&=& O(e^x) \cr
\sin(x) &\ll& 1 \cr
(x+2)^{10}&\ll& x^{10}\cr
(x+2)^{10}&=& x^{10} + O(x^9) . \cr
\end{eqnarray}

Examples assuming $x\to 0$:
\begin{eqnarray}
x^4 &\ll& x^3 \cr
\log(1+x)&=& O(x)   \cr
\log(1+x)&=& x +O(x^2)  \cr
\sin(x) &\ll& x \cr
(x+2)^{10}&=& 1024+O(x)  . \cr
\end{eqnarray}

\subsection{Little-$o$ notation}

The statement
\begin{equation}
f(x)=o(g(x))
\ \ \ \ \ \ \ \ \ \ \ \ 
\text{as }\ \ x\to\infty
\end{equation}
is pronounced ``$f(x)$ is little oh of $g(x)$.''  It means 
\begin{equation}
\lim_{x\to\infty} \frac{f(x)}{g(x)} =0.
\end{equation}
Equivalently, for all $C>0$, if $x$ is sufficiently large then
$|f(x)|\le C g(x)$.  It is like big-$O$ where the implied constant
can be made arbitrarily small.   

Note that $f(x)\sim g(x)$,
``$f(x)$ is asymptotic to $g(x)$'' is equivalent to
$f(x)=(1+o(1)) g(x)$.


\begin{thebibliography}{99}

\bibitem[CF]{CF}
\newblock
J. B. Conrey and D. W. Farmer,  \emph{Mean values of $L$-functions and symmetry}, 
\jour  Internat. Math. Res. Notices
\yr 2000 {\bf 17}  pp. 883--908.

\bibitem[CG]{CG}
J. B. Conrey and A. Ghosh,
 {\it  Mean values of the Riemann zeta-function},
\jour Mathematika
\vol 31
\yr 1984
\pages  159--161

\bibitem[GG]{GG}
\newblock
D.A. Goldston and S. M. Gonek, \emph{Mean value theorems for long 
Dirichlet polynomials and tails of Dirichlet series},  Acta Arith.  84  
(1998),  no. 2, 155--192. 

\bibitem[KS]{KS}
\newblock
 J. P. Keating and N. C. Snaith,  \emph{Random matrix theory and $L$-functions at $s=\frac12 $}, \jour Comm. Math. Phys.
 {\bf 214} \yr 2000
 pp.  91--110.

\bibitem[MV]{MV}
H. L. Montgomery and R. C. Vaughan
\emph{The large sieve},
\jour Mathematika, 1973 {\bf 20}, 119--134 

\bibitem[T]{T}
\newblock
 E. C. Titchmarsh,  The theory of the Riemann zeta-function. Second edition. Edited and with a preface by D. R. Heath-Brown, \publ The Clarendon Press, Oxford University Press
\publaddr New York
\yr 1986.

\end{thebibliography}
\end{document}